\documentclass[11pt,a4paper]{amsart}

\usepackage{amsfonts}

\usepackage[T1]{fontenc}

\usepackage{enumerate}

\usepackage[abbrev,backrefs]{amsrefs}

\usepackage[utf8]{inputenc}

\usepackage{amssymb}

\usepackage{palatino}

\usepackage[nameinlink]{cleveref}

\newcommand{\abs}[1]{\mathopen\lvert#1\mathclose\rvert}

\newcommand{\norm}[1]{\mathopen\lVert#1\mathclose\rVert}

\newcommand{\N}{{\mathbb N}}
\newcommand{\R}{{\mathbb R}}

\newcommand{\cM}{\mathcal{M}}

\DeclareMathOperator{\Div}{div}

\newcommand{\dif}{\,\mathrm{d}}

\newcommand{\meas}[1]{\left| #1 \right|}

\theoremstyle{plain}

\newtheorem{theorem}{Theorem}

\theoremstyle{definition}

\theoremstyle{remark}

\newtheorem{example}{Example}

\numberwithin{equation}{section}

\title{Flat solutions of the 1-Laplacian equation}

\author{Luigi Orsina}
\address{
Luigi Orsina\hfill\break\indent
``Sapienza'' Universit\`a di Roma\hfill\break\indent
Dipartimento di Matematica \hfill\break\indent
P.le A.~Moro 2\hfill\break\indent
00185 Roma, Italy}

\author{Augusto C. Ponce}
\address{
Augusto C. Ponce\hfill\break\indent
 Universit{\'e} catholique de Louvain\hfill\break\indent
 Institut de Recherche en Math{\'e}matique et Physique\hfill\break\indent
 Chemin du cyclotron 2, bte L7.01.02\hfill\break\indent
1348 Louvain-la-Neuve, Belgium}



\begin{document}

\begin{abstract}
For every \(f \in L^N(\Omega)\) defined in an open bounded subset \(\Omega\) of \(\R^N\), we prove that a solution \(u \in W_0^{1, 1}(\Omega)\) of the \(1\)-Laplacian equation \({-}\Div{\big(\frac{\nabla u}{|\nabla u|}\big)} = f\) in \(\Omega\) satisfies \(\nabla u = 0\) on a set of positive Lebesgue measure.
The same property holds if \(f \not\in L^N(\Omega)\) has small norm in the Marcinkiewicz space of weak-\(L^{N}\) functions or if \(u\) is a BV minimizer of the associated energy functional.
The proofs rely on Stampacchia's truncation method.
\end{abstract}

\dedicatory{À Jean Mawhin, dont l'enthousiasme et l'amour pour les mathématiques font encore rêver les nouvelles générations.}

\maketitle

\section{Introduction}

Let \(\Omega \subset \R^N\) be a bounded smooth open subset.
Given a convex function \(\Phi : \R^{N} \to \R\) and \(f \in L^{1}(\Omega)\), consider the energy functional
\[{}
E_{\Phi}(u)
= \int\limits_{\Omega} \Phi(\nabla u) - \int\limits_{\Omega} f u,
\]
defined on some class of functions \(u : \Omega \to \R\) for which the integrands are summable.
Although \(\Phi\) need not be smooth, one can express the Euler--Lagrange equation of \(E_{\Phi}\) using the subdifferential of \(\Phi\).{}
Indeed, by convexity of \(\Phi\), at each point \(x \in \R^{N}\) there exists a subgradient \(\xi \in \R^{N}\) such that
\[{}
\Phi(y) \ge \Phi(x) + \xi \cdot (y - x),
\]
for every \(y \in \R^{N}\); see \cite{Mawhin-Willem}*{Chapter~2}.
Denoting the collection of all subgradients \(\xi\) at \(x\) by \(\partial \Phi(x)\), one can then write the Euler--Lagrange equation of \(E_{\Phi}\) at some function \(u\) as~\citelist{\cite{Ekeland-Temam}*{Chapter~IV} \cite{Papageorgiou}}
\begin{equation}
\label{eqEulerLagrange1}
- \Div{Z} = f	\quad \text{in the sense of distributions in \(\Omega\),}
\end{equation}
where \(Z\) is a summable function with values in \(\R^{N}\) such that 
\begin{equation}
\label{eqEulerLagrange2}
Z \in \partial\Phi(\nabla u) \quad \text{almost everywhere in \(\Omega\).}
\end{equation}

For example, if \(\Phi_{p}(x) = \abs{x}^{p}/p\) for some exponent \(p > 1\), then \(\Phi_{p}\) is differentiable pointwise.
Thus, \(\partial\Phi_{p}(x) = \{\abs{x}^{p-2}x\}\), and one recovers an equation involving the \(p\)-Laplace operator:
\[{}
-\Delta_{p}{u} = -\Div{(\abs{\nabla u}^{p-2}\nabla u)} = f.
\]
When \(p = 1\), the function \(\Phi_{1}\) is not differentiable at \(0\), and one should be careful about the meaning of the quotient \(\nabla u/\abs{\nabla u}\) that appears in the formal notation of the \(1\)-Laplacian.{}
The correct interpretation is based on the formalism of subdifferentials above.
Indeed, for \(\Phi_{1}(x) = \abs{x}\), one has
\begin{equation}
	\label{eqSubdifferential1}
\partial\Phi_{1}(x){}
=
\begin{cases}
	\overline{B_{1}(0)}		& \text{if \(x = 0\),}\\
	\left\{x/\abs{x}\right\}	& \text{if \(x \ne 0\),}
\end{cases}
\end{equation}
where \(B_1(0)\) denotes the unit open ball in \(\R^N\).{}
The vector field \(Z\) in the Euler--Lagrange equation now satisfies the conditions:
\[{}
\abs{Z} \le 1
\quad \text{and} \quad 
Z \abs{\nabla u} = \nabla u
\]
almost everywhere in \(\Omega\).{}
Observe that, in dimension \(1\), equation~\eqref{eqSubdifferential1} provides one with the maximal monotone graph associated to the sign function.

Assuming that \(f \in L^{N}(\Omega)\), the functional \(E_{\Phi_{1}}\) associated to \(\Phi_{1}\) is well-defined in \(W_{0}^{1, 1}(\Omega)\), and the Euler--Lagrange equation \eqref{eqEulerLagrange1}--\eqref{eqEulerLagrange2} is indeed satisfied by a minimizer.
The goal of this paper is to show that one cannot abandon the vector field \(Z\) and replace it by the quotient \(\nabla u/\abs{\nabla u}\) since the gradient \(\nabla u\) must vanish on a set of positive Lebesgue measure.

Every function \(u \in W^{1, 1}(\Omega)\) such that \(\nabla u \ne 0\) a.e.~in \(\Omega\) has a legitimate 1-Laplacian \(\Delta_1 u\)
defined in the sense of distributions as
\[
\langle  \Delta_1 u, \varphi \rangle := - \int\limits_\Omega \frac{\nabla u }{|\nabla u|} \cdot \nabla \varphi,
\]
for every test function \(\varphi \in C_c^\infty(\Omega)\) with compact support in \(\Omega\), but even for smooth functions \(u\) something strange happens near an interior extremum point:

\begin{example}
\label{exampleSmooth}
Let \(B_1(0)\) be the unit ball in \(\R^N\) and let \(u : B_1(0) \to \R\) be the function defined by \(u(x) = 1 - |x|^2\).
In the sense of distributions we have, for \(N = 1\), 
\[
-\Delta_{1} u = 2 \delta_{0},
\]
while for \(N \ge 2\),
\[
-\Delta_{1} u = \frac{N-1}{|x|}.
\]
\end{example}

In the previous example, the topological singularity of the vector field \(-x/\abs{x}\) is detected by its divergence, and
 the \(1\)-Laplacian does not belong to \(L^N(\Omega)\).
We show that this is a general fact that holds for Sobolev functions, not necessarily smooth:

\begin{theorem}
\label{theoremMain}
There exists no function \(u \in W_0^{1, 1}(\Omega)\) such that \(\nabla u \ne 0\) a.e.~in \(\Omega\) and 
\[
\Delta_{1} u  \in L^N(\Omega).
\]
\end{theorem}

In \Cref{exampleSmooth} above for \(N \ge 2\), one sees that the right-hand side belongs to the Marcinkiewicz space \(\cM^{N}(\Omega)\) of weak-\(L^{N}\) functions \(f\) in \(\Omega\) equipped with the seminorm
\[{}
\norm{f}_{\cM^{N}(\Omega)}
= \sup_{A \subset \Omega}{\frac{1}{\meas{A}^{\frac{N-1}{N}}} \int_{A}{\abs{f}}},
\]
where \(\meas{A}\) denotes the Lebesgue measure of \(A\) and the supremum is computed with respect to every Borel subset of \(\Omega\).
In the case of the example, the function \(f = (N-1)/\abs{x}\) satisfies
\begin{equation}
\label{eqMarcinkiewiczCritical}
\norm{f}_{\cM^{N}(B(0; 1))} 
= N \omega_{N}^{1/N},
\end{equation}
where \(\omega_{N}\) denotes the volume of the unit ball in \(\R^{N}\).

A variant of the proof of Theorem~\ref{theoremMain} based on Peetre--Alvino's imbedding of \(W^{1, 1}(\R^{N})\) in the Lorentz space \(L^{\frac{N}{N-1}, 1}(\R^{N})\) shows that this quantity \eqref{eqMarcinkiewiczCritical} is critical for the existence of flat levels of solutions involving the \(1\)-Laplacian:

\begin{theorem}
\label{theoremMainLorentz}
Let \(N \ge 2\).{}
There exists no function \(u \in W_0^{1, 1}(\Omega)\) such that \(\nabla u \ne 0\) a.e.~in \(\Omega\),
\[
\Delta_{1} u \in \cM^N(\Omega) \quad \text{and} \quad \norm{\Delta_{1} u}_{\cM^{N}(\Omega)} < N \omega_{N}^{1/N}.
\]
\end{theorem}

\Cref{theoremMain,theoremMainLorentz} are related to the degenerate limit behavior of solutions of the \(p\)-Laplacian equation as \(p\) tends to~\(1\) that has been studied by several authors~\cites{Cicalese-Trombetti,Mercaldo-Segura-Trombetti-2008,Mercaldo-Segura-Trombetti-2009}, starting with the pioneering work of Kawohl~\cites{Kawohl}, and also clarify the need for relying on the vector field \(Z\) in replacement of \(\nabla u/\abs{\nabla u}\).{}

\begin{example}
For any \(0 < r < 1\), let \(u : B_1(0) \to \R\) be the function defined by
\[
u(x) =
\begin{cases}
1 - |x|^2	& \text{if \(|x| \ge r\),}\\
1 - r^2	& \text{if \(|x| < r\).}
\end{cases}
\]
Then, \(u \in W_0^{1, 1}(B_1(0))\).
If \(Z : B_1(0) \to \R^N\) is any smooth extension of the function 
\[{}
x \in B_1(0) \setminus B_r(0) \longmapsto - \frac{x}{|x|} \in \R^N,{}
\]
then \(u\) and \(Z\) satisfy the Euler--Lagrange equation \eqref{eqEulerLagrange1}--\eqref{eqEulerLagrange2} for some \(f \in L^\infty(B_1(0))\).
\end{example}

Observe that the Sobolev space \(W_{0}^{1, 1}(\Omega)\) is not the natural setting for looking for minimizers of \(E_{\Phi_{1}}\),
due to the lack of reflexivity of \(L^{1}(\Omega; \R^{N})\).{}
This is in contrast with minimization problems in \(W^{1, p}(\Omega)\) for \(1 < p < +\infty\) which can be investigated using techniques based on the uniform convexity of the space; see \cite{Dinca-Jebelean-Mawhin}.

Let us assume that \(E_{\Phi_{1}}\) is bounded from below for some given \(f \in L^{N}(\Omega)\).{}
This is the case for example if the norm \(\norm{f}_{L^{N}(\Omega)}\) is small, depending on the Sobolev constant; see e.g.~\cite{Kawohl-Schuricht}.
One can now take a minimizing sequence \((u_{n})_{n \in \N}\) in \(W_{0}^{1, 1}(\Omega)\) such that
\[{}
\lim_{n \to \infty}{E_{\Phi_{1}}(u_{n})}
= \inf_{W_{0}^{1, 1}(\Omega)}{E_{\Phi_{1}}}.
\] 
Each function \(u_{n}\), extended by zero to \(\R^{N}\), is an element of \(W^{1, 1}(\R^{N})\).{}
Since the sequence \((\nabla u_{n})_{n \in \N}\) is bounded in \(L^{1}(\R^{N}; \R^{N})\), we can extract a subsequence \((\nabla u_{n_{k}})_{k \in \N}\) converging weakly to some finite vector-valued measure in \(\R^{N}\) supported in \(\overline{\Omega}\).{}
Applying the Rellich--Kondrashov compactness theorem, we deduce that there exists \(u \in BV(\R^{N})\) such that \(u = 0\) in \(\R^{N} \setminus \Omega\), and
\[{}
\lim_{k \to \infty}{E_{\Phi_{1}}(u_{n_{k}})}
\ge \int_{\R^{N}} \abs{Du} - \int_{\Omega} fu.
\]

The limit function \(u\) is a minimizer of the extended functional
\begin{equation}
\label{eqFunctionalExtended}
\overline E_{\Phi_{1}}(v)
:= \int_{\R^{N}} \abs{Dv} - \int_{\Omega} fv,
\end{equation}
over the class of functions \(v \in BV(\R^{N})\) such that \(v = 0\) in \(\R^{N} \setminus \Omega\).{}
Such a functional provides a relaxed formulation of the minimization problem for which a solution exists; see~\cite{Giaquinta-Modica-Soucek}.
In the spirit of \Cref{theoremMain,theoremMainLorentz}, minimizers of \eqref{eqFunctionalExtended} must have flat level sets:

\begin{theorem}
	\label{theoremMainBV}
	Let \(f \in L^{N}(\Omega)\) and let \(u \in BV(\R^{N})\) with \(u = 0\) in \(\R^{N} \setminus \Omega\) be a minimizer of the extended functional \(\overline E_{\Phi_{1}}\).{}
	Then, \(u \in L^{\infty}(\R^{N})\) and the set of extremal points 
	\[{}
	\big\{x \in \R^{N} : \abs{u(x)} = \norm{u}_{L^{\infty}(\R^{N})}\big\}
	\]
	has positive Lebesgue measure.
\end{theorem}

We deduce in this case that the absolute continuous part \(D^{\mathrm{a}}u\) of the measure \(Du\) vanishes a.e.~on a set of positive measure, since \(D^{\mathrm{a}}u = 0\) a.e.~on level sets \(\{u = \alpha\}\) for every \(\alpha \in \R\) \cite{Ambrosio-Fusco-Pallara}*{Proposition~3.73}.
The counterpart of Theorem~\ref{theoremMainBV} involving the condition \(\norm{f}_{\cM^{N}(\Omega)} < N \omega_{N}^{1/N}\) is true but uninteresting since \(\overline E_{\Phi_{1}}\) is nonnegative and \(0\) is the unique minimizer.
This follows from a standard application of Alvino's version of the Sobolev inequality in Lorentz spaces.

Renormalized solutions to equations involving the \(1\)-Laplacian have been introduced in the spirit of the relaxed minimization problem above, but in general such solutions merely have bounded variation or do not satisfy the homogeneous Dirichlet boundary condition \cites{Mazon-Segura,Andreu-Ballester-Caselles-Mazon,Andreu-DallAglio-Segura,Demengel-2004,Bellettini-Caselles-Novaga,Alter-Caselles-Chambolle}.

\begin{example}[\cite{Mazon-Segura}*{Remark~3.10}]
For every \(N < r \le R\), the function \(u = (1 - N/r)\chi_{B_r(0)}\) is a renormalized solution of the Dirichlet problem
\[
\left\{
\begin{alignedat}{2}
- \Delta_1 v & = h - v	\quad && \text{in \(B_R(0)\),}\\
v & = 0					\quad && \text{on \(\partial B_R(0)\),}
\end{alignedat}
\right.
\]
with bounded datum \(h = \chi_{B_r(0)}\). 
Note that if \(r < R\) then \(u_{r}\) is a BV function with compact support in \(B_{R}(0)\), while if \(r = R\) then \(u_{r}\) is a \(W^{1,1}\) function which does not vanish on the boundary.
\end{example}

In the next section, we prove \Cref{theoremMain,theoremMainLorentz,theoremMainBV}.
This paper is a revised and extended version of a note written by the authors in 2012 that was only available at the \texttt{arxiv.org} website.

\section{Proofs of the main results}

\begin{proof}[Proof of \Cref{theoremMain}]
Assume by contradiction that there exists a function \(u \in W_{0}^{1, 1}(\Omega)\) such that \(\nabla u \ne 0\) almost everywhere in \(\Omega\) and \(f := \Delta_{1} u \in L^N(\Omega)\).{}
Then,
\[
\int\limits_\Omega \frac{\nabla u}{|\nabla u|} \cdot \nabla \varphi
=
\int\limits_\Omega f \varphi,
\]
for every \(\varphi \in C_c^\infty(\Omega)\).
Note that \({\nabla u}/{|\nabla u|} \in L^\infty(\Omega)\) and \(u \in L^{\frac{N}{N-1}}(\Omega)\) by the Gagliardo-Nirenberg-Sobolev imbedding.
By density of \(C_c^\infty(\Omega)\) in \(W_0^{1, 1}(\Omega)\) we deduce that
\begin{equation}\label{eq:TestFunctionW11}
\int\limits_\Omega \frac{\nabla u}{|\nabla u|} \cdot \nabla v
=
\int\limits_\Omega f v,
\end{equation}
for every \(v \in W_0^{1, 1}(\Omega)\).

We proceed using Stampacchia's truncation method.
For this purpose, for every \(\kappa > 0\) let \(G_\kappa : \R \to \R\) be the function defined by
\begin{equation}
\label{eqTruncationStampacchia}
G_\kappa(t) =
\begin{cases}
t + \kappa	& \text{if \(t < -\kappa\),}\\
0	& \text{if \(- \kappa \le t \le \kappa\),}\\
t - \kappa	& \text{if \(t > \kappa\).}
\end{cases}
\end{equation}
Since \(u \in W_0^{1, 1}(\Omega)\), we have \(G_\kappa(u) \in W_0^{1, 1}(\Omega)\).
Hence,
\[
\frac{\nabla u}{|\nabla u|} \cdot \nabla G_\kappa(u) 
= G_\kappa'(u) |\nabla u|
= |\nabla G_\kappa(u)| .
\]
Applying identity~\eqref{eq:TestFunctionW11} with test function \(G_\kappa(u)\), we get
\[
\int\limits_\Omega |\nabla G_\kappa(u)|
=
\int\limits_\Omega f G_\kappa(u).
\]
Since \(G_\kappa\) vanishes on the interval \([-\kappa, \kappa]\), by the H\"older inequality we have
\[
\int\limits_\Omega f G_\kappa(u)
=
\int\limits_{\{\abs{u} > \kappa\}} f G_\kappa(u)
\le
\|f\|_{L^N(\{\abs{u} > \kappa\})} \|G_\kappa(u)\|_{L^{\frac{N}{N-1}}(\Omega)}.
\]
Thus,
\[
\int\limits_\Omega |\nabla G_\kappa(u)|
\le
\|f\|_{L^N(\{\abs{u} > \kappa\})} \|G_\kappa(u)\|_{L^{\frac{N}{N-1}}(\Omega)}.
\]
By the Gagliardo-Nirenberg-Sobolev inequality, 
\[
\|G_\kappa(u)\|_{L^\frac{N}{N-1}(\Omega)} 
\le
C \int\limits_\Omega |\nabla G_\kappa(u)|,
\]
for some constant \(C > 0\) depending only on the dimension \(N\).
Hence,
\begin{equation}
\label{eq:inequalityLN}
\big(1 - C\|f\|_{L^N(\{\abs{u} > \kappa\})}\big) \|G_\kappa(u)\|_{{L^{\frac{N}{N-1}}(\Omega)}} \le 0.
\end{equation}

Let \(T := \|u\|_{L^\infty(\Omega)}\) if \(u\) is essentially bounded, or \(T := +\infty\) otherwise.
We have
\[
\lim_{\kappa \nearrow T}{\|f\|_{L^N(\{\abs{u} > \kappa\})}} 
= \|f\|_{L^N(\{\abs{u} = T\})}.
\]
We observe that the set \(\{\abs{u} = T\}\) has zero Lebesgue measure. 
This is indeed the case when \(T = +\infty\) since \(u\) is finite a.e.
When \(T < +\infty\), we observe that \(\nabla u = 0\) a.e.~on the level set \(\{\abs{u} = T\}\); since by assumption \(\nabla u \ne 0\) a.e.~in \(\Omega\), the set \(\{u = T\}\) must have zero Lebesgue measure.
This implies that
\[
\lim_{\kappa \nearrow T}{\|f\|_{L^N(\{\abs{u} > \kappa\})}} 
= \|f\|_{L^N(\{\abs{u} = T\})} = 0.
\]
In particular, there exists \(0 < \kappa < T\) such that \(C \|f\|_{L^N(\{\abs{u} > \kappa\})} < 1\).
We deduce from \eqref{eq:inequalityLN} that 
\[
\|G_\kappa(u)\|_{{L^{\frac{N}{N-1}}(\Omega)}} \le 0.
\]
Therefore, \(\abs{u} \le \kappa\) a.e.~in \(\Omega\).
Hence, \(T = \norm{u}_{L^\infty(\Omega)} \le \kappa\), and this contradicts the choice of \(\kappa\).
The proof of the theorem is complete.
\end{proof}

To prove Theorem~\ref{theoremMainLorentz}, we rely on Peetre's imbedding of Sobolev functions in Lorentz spaces, with the best constant computed by Alvino.
We recall that the Lorentz space \(L^{p, 1}(\R^{N})\) for \(1 \le p < \infty\) can be defined as the vector space of measurable functions \(g\) in \(\R^{N}\) such that
\[{}
\norm{g}_{L^{p, 1}(\R^{N})}
:= \int_{0}^{\infty} \meas{\{\abs{g} > t\}}^{1/p} \dif t < +\infty.
\]
Using an equivalent definition to this one, Lorentz~\cite{Lorentz} established the duality between \(L^{p, 1}(\R^{N})\) and \(\cM^{\frac{p}{p-1}}(\R^{N})\) for \(p > 1\) by proving an estimate which amounts to
\[{}
\int_{\R^{d}}{\abs{fg}}
\le \norm{f}_{\cM^{\frac{p}{p-1}}(\R^{N})} \norm{g}_{L^{p, 1}(\R^{N})} ,
\]  
for every \(g \in L^{p, 1}(\R^{N})\) and \(f \in \cM^{\frac{p}{p-1}}(\R^{N})\), where
\[{}
\norm{f}_{\cM^{\frac{p}{p-1}}(\R^{N})}
:= \sup_{A \subset \Omega}{\frac{1}{\meas{A}^{\frac{1}{p}}} \int_{A}{\abs{f}}}\ ;
\]
see \cite[Theorem~5]{Lorentz} and the computation of the Lorentz norm in \cite[{Section~2}]{Bahouri-Cohen}.
Here one should not rely on the quasi-norm \(\sup_{t > 0}{\bigl\{t \meas{\{\abs{f} > t\}}^{\frac{p-1}{p}} \bigr\}}\), which gives a quantity that is only equivalent to \(\norm{f}_{\cM^{\frac{p}{p-1}}(\R^{N})}\).

Peetre~\cite{Peetre} proved by interpolation that \(W^{1, 1}(\R^{N}) \subset L^{\frac{N}{N-1}, 1}(\R^{N})\) and Alvino~\cite{Alvino:1977} later showed using rearrangements that the inequality
\[{}
\norm{v}_{L^{\frac{N}{N-1}, 1}(\R^{N})}
\le \gamma_{1} \norm{\nabla v}_{L^{1}(\R^{N})}
\]
holds with the best constant given by \(\gamma_{1} := 1/(N\omega_{N}^{1/N})\).{}

\begin{proof}[Proof of Theorem~\ref{theoremMainLorentz}]
	Proceeding as in the previous proof, by the duality between \(L^{\frac{N}{N-1}, 1}\) and \(\cM^{N}\) one gets
\[
\int\limits_{\R^{N}} |\nabla G_\kappa(u)|{}
= \int_{\Omega} f G_{\kappa}(u)
\le
\|f\|_{\cM^{N}(\R^{N})} \|G_\kappa(u)\|_{L^{\frac{N}{N-1}, 1}(\R^{N})},
\]
where the functions \(f\) and \(u\) have been extended by zero to \(\R^{N}\); this does not change their seminorms.
Using Alvino's improvement of the Sobolev inequality with \(v = G_{\kappa}(u)\), it follows that
\[{}
\big(1 - \gamma_{1}\|f\|_{\cM^N(\R^{N})}\big) \|G_\kappa(u)\|_{{L^{\frac{N}{N-1}, 1}(\R^{N})}} \le 0.
\]
Under the assumption of the theorem we have \(\|f\|_{\cM^N(\R^{N})} < 1/\gamma_{1}\), hence the quantity in parenthesis is positive.
We deduce that \(\|G_\kappa(u)\|_{{L^{\frac{N}{N-1}, 1}(\R^{N})}} = 0\) for every \(\kappa > 0\), and then \(u = 0\) a.e.~in \(\Omega\), but this is not possible.
\end{proof}

The proof of \Cref{theoremMainBV} relies on a property of BV~function related to the chain rule.
For this purpose, given \(\kappa > 0\) denote by \(T_{\kappa} : \R \to \R\) the truncation function at levels \(\pm \kappa\):{}
\[
T_\kappa(t) =
\begin{cases}
- \kappa	& \text{if \(t < -\kappa\),}\\
t	& \text{if \(- \kappa \le t \le \kappa\),}\\
\kappa	& \text{if \(t > \kappa\).}
\end{cases}
\]
Observe that, for every \(t \in \R\),{}
\begin{equation}
\label{eqIdentity}
t = T_{\kappa}(t) + G_{\kappa}(t),
\end{equation}
where \(G_{\kappa}\) is the function defined by \eqref{eqTruncationStampacchia}.
Since \(T_{\kappa}\) and \(G_{\kappa}\) are Lipschitz continuous, it is straightforward to verify using an approximation argument that \(T_{\kappa}(u)\) and \(G_{\kappa}(u)\) both belong to \(BV(\R^{N})\) for every \(u \in BV(\R^{N})\).{}
In addition, by the identity above we have
\[{}
Du = D(T_{\kappa}(u)) + D(G_{\kappa}(u)).
\]

One then verifies that
\begin{equation}
\label{eqIdentityMeasures}
\int_{\R^{N}}{\abs{Du}} 
= \int_{\R^{N}}{\abs{D(T_{\kappa}(u))}} + \int_{\R^{N}}{\abs{D(G_{\kappa}(u))}},
\end{equation}
where, for a given vector-valued measure \(\mu\),{}
\[{}
\int_{\R^{N}}{\abs{\mu}}
= \sup{\bigg\{ \int_{\R^{N}} \Phi \cdot \mu : \Phi \in C_{c}^{\infty}(\R^{N}; \R^{N})\ \text{and}\ \abs{\Phi} \le 1\ \text{in \(\R^{N}\)}  \bigg\}}.
\]

The inequality \(\le\) in \eqref{eqIdentityMeasures} follows from the triangle inequality in \(\R^{N}\).
The reverse inequality \(\ge\) can be deduced from Vol'pert's chain rule for \(BV\) functions~\cite{Ambrosio-DalMaso}.
A more elementary approach is based on an approximation of \(u\) using the sequence of smooth functions \((\rho_{n} * u)_{n \in \N}\), where \((\rho_{n})_{n \in \N}\) is a sequence of mollifiers in \(C_{c}^{\infty}(\R^{N})\).{}
In this case, one observes that
\[{}
\int_{\R^{N}}{\abs{D(\rho_{n} * u)}} \to \int_{\R^{N}}{\abs{Du}}
\]
as \(n \to \infty\); see \cite{Evans-Gariepy}*{Theorem~5.3}.
On the other hand, there exist a subsequence \((\rho_{n_{j}} * u)_{j \in \N}\) and finite positive measures \(\sigma_{1}\) and \(\sigma_{2}\) such that
\begin{alignat*}{2}
\abs{D(T_{\kappa}(\rho_{n_{j}} * u))} 
& \overset{*}{\rightharpoonup} \sigma_{1} && \quad \text{in \(\cM(\R^{N}; \R^{N})\),}\\
\abs{D(G_{\kappa}(\rho_{n_{j}} * u))} 
& \overset{*}{\rightharpoonup} \sigma_{2} && \quad \text{in \(\cM(\R^{N}; \R^{N})\),}
\end{alignat*}
as \(j \to \infty\), where \(\sigma_{1} \ge \abs{D(T_{\kappa}(u))}\) and \(\sigma_{2} \ge \abs{D(G_{\kappa}(u))}\).{}
This implies the reverse inequality in \eqref{eqIdentityMeasures}.

\begin{proof}[Proof of \Cref{theoremMainBV}]
	Since \(u\) minimizes \(\overline E_{\Phi_{1}}\), and \(T_{\kappa}(u)\) is also an admissible function in the minimization class, we have
	\[{}
	\overline E_{\Phi_{1}} (u)
	\le \overline E_{\Phi_{1}}(T_{\kappa}(u)).
	\]
	Thus,
	\[{}
	\int_{\R^{N}}{\big[\abs{Du} - \abs{D(T_{\kappa}(u))}\big]}
	\le \int_{\R^{N}}{f(u - T_{\kappa}(u))}.
	\]
	We deduce from \eqref{eqIdentityMeasures} and \eqref{eqIdentity} that
	\[{}
	\int_{\R^{N}}{\abs{D(G_{\kappa}(u))}}
	\le \int_{\R^{N}}{f G_{\kappa}(u)}.
	\]
	We can now pursue the strategy of the proof of \Cref{theoremMain} to get the conclusion.
	Indeed, the Sobolev and Hölder inequalities imply that
	\[{}
	\big(1 - C\|f\|_{L^N(\{\abs{u} > \kappa\})}\big) \|G_\kappa(u)\|_{{L^{\frac{N}{N-1}}(\Omega)}} \le 0.	
	\]
	For every \(0 < \kappa < \norm{u}_{L^{\infty}(\R^{N})}\), we have \(\|G_\kappa(u)\|_{{L^{\frac{N}{N-1}}(\Omega)}} > 0\).{}
	Thus,
	\[{}
	\|f\|_{L^N(\{\abs{u} > \kappa\})} \ge \frac{1}{C}.
	\]
	Since \(u\) is finite a.e., this inequality cannot hold for every \(\kappa > 0\).{}
	Therefore, \(\norm{u}_{L^{\infty}(\R^{N})} < \infty\).{}
	Letting \(\kappa \to \norm{u}_{L^{\infty}(\R^{N})}\), we deduce that \(\{\abs{u} > \kappa\}\) has positive measure.
\end{proof}




\section*{Acknowledgements}
The authors would like to thank Francesco Petitta for bringing to their attention problems related to the \(1\)-Laplacian.
The second author (ACP) was supported by the Fonds de la Recherche scientifique--FNRS under research grants 1.5.199.10F and J.0026.15. 
He warmly thanks the Dipartimento di Matematica of ``Sapienza'' Universit\`a di Roma for the invitation and hospitality.



\begin{bibdiv}

\begin{biblist}


\bib{Alter-Caselles-Chambolle}{article}{
   author={Alter, F.},
   author={Caselles, V.},
   author={Chambolle, A.},
   title={A characterization of convex calibrable sets in $\Bbb R^N$},
   journal={Math. Ann.},
   volume={332},
   date={2005},
   pages={329--366},
}
			
	
\bib{Alvino:1977}{article}{
   author={Alvino, Angelo},
   title={Sulla diseguaglianza di Sobolev in spazi di Lorentz},
   journal={Boll. Un. Mat. Ital. A (5)},
   volume={14},
   date={1977},
   pages={148--156},
}

\bib{Ambrosio-DalMaso}{article}{
   author={Ambrosio, L.},
   author={Dal Maso, G.},
   title={A general chain rule for distributional derivatives},
   journal={Proc. Amer. Math. Soc.},
   volume={108},
   date={1990},
   pages={691--702},
}

\bib{Ambrosio-Fusco-Pallara}{book}{
   author={Ambrosio, Luigi},
   author={Fusco, Nicola},
   author={Pallara, Diego},
   title={Functions of bounded variation and free discontinuity problems},
   series={Oxford Mathematical Monographs},
   publisher={The Clarendon Press, Oxford University Press, New York},
   date={2000},
}

\bib{Andreu-DallAglio-Segura}{article}{
   author={Andreu, Fuensanta},
   author={Dall'Aglio, Andrea},
   author={Segura de Le{\'o}n, Sergio},
   title={Bounded solutions to the 1-Laplacian equation with a critical
   gradient term},
   journal={Asymptot. Anal.},
   volume={80},
   date={2012},
   pages={21--43},
}

\bib{Andreu-Ballester-Caselles-Mazon}{article}{
   author={Andreu, F.},
   author={Ballester, C.},
   author={Caselles, V.},
   author={Maz{\'o}n, J. M.},
   title={The Dirichlet problem for the total variation flow},
   journal={J. Funct. Anal.},
   volume={180},
   date={2001},
   pages={347--403},
}


\bib{Bahouri-Cohen}{article}{
   author={Bahouri, Hajer},
   author={Cohen, Albert},
   title={Refined Sobolev inequalities in Lorentz spaces},
   journal={J. Fourier Anal. Appl.},
   volume={17},
   date={2011},
   pages={662--673},
}


\bib{Bellettini-Caselles-Novaga}{article}{
   author={Bellettini, Giovanni},
   author={Caselles, Vicent},
   author={Novaga, Matteo},
   title={Explicit solutions of the eigenvalue problem $-{\rm
   div}\big(\frac{Du}{\vert Du\vert }\big)=u$ in $\mathbb{R}^2$},
   journal={SIAM J. Math. Anal.},
   volume={36},
   date={2005},
   pages={1095--1129},
}


\bib{Cicalese-Trombetti}{article}{
   author={Cicalese, Marco},
   author={Trombetti, Cristina},
   title={Asymptotic behaviour of solutions to $p$-Laplacian equation},
   journal={Asymptot. Anal.},
   volume={35},
   date={2003},
   pages={27--40},
}


\bib{Demengel-2004}{article}{
   author={Demengel, Fran{\c{c}}oise},
   title={Functions locally almost 1-harmonic},
   journal={Appl. Anal.},
   volume={83},
   date={2004},
   pages={865--896},
}


\bib{Dinca-Jebelean-Mawhin}{article}{
   author={Dinca, G.},
   author={Jebelean, P.},
   author={Mawhin, J.},
   title={Variational and topological methods for Dirichlet problems with
   $p$-Laplacian},
   journal={Port. Math. (N.S.)},
   volume={58},
   date={2001},
   pages={339--378},
}

\bib{Ekeland-Temam}{book}{
   author={Ekeland, Ivar},
   author={T{\'e}mam, Roger},
   title={Convex analysis and variational problems},
   series={Classics in Applied Mathematics},
   volume={28},
   publisher={Society for Industrial and Applied Mathematics (SIAM),
   Philadelphia, PA},
   date={1999},
}


\bib{Evans-Gariepy}{book}{
   author={Evans, Lawrence C.},
   author={Gariepy, Ronald F.},
   title={Measure theory and fine properties of functions},
   series={Textbooks in Mathematics},
   edition={Revised edition},
   publisher={CRC Press},
   address={Boca Raton, FL},
   date={2015},
}

\bib{Giaquinta-Modica-Soucek}{article}{
   author={Giaquinta, M.},
   author={Modica, G.},
   author={Sou{\v{c}}ek, J.},
   title={Functionals with linear growth in the calculus of variations. I,
   II},
   journal={Comment. Math. Univ. Carolin.},
   volume={20},
   date={1979},
   pages={143--156, 157--172},
}



\bib{Kawohl}{article}{
   author={Kawohl, Bernhard},
   title={On a family of torsional creep problems},
   journal={J. Reine Angew. Math.},
   volume={410},
   date={1990},
   pages={1--22},
}

\bib{Kawohl-Schuricht}{article}{
   author={Kawohl, Bernd},
   author={Schuricht, Friedemann},
   title={Dirichlet problems for the 1-Laplace operator, including the
   eigenvalue problem},
   journal={Commun. Contemp. Math.},
   volume={9},
   date={2007},
   pages={515--543},
}


\bib{Lorentz}{article}{
   author={Lorentz, G. G.},
   title={Some new functional spaces},
   journal={Ann. of Math. (2)},
   volume={51},
   date={1950},
   pages={37--55},
}

\bib{Mazon-Segura}{article}{
   author={Maz{\'o}n, Jos{\'e} M.},
   author={Segura de Le{\'o}n, Sergio},
   title={The Dirichlet problem for a singular elliptic equation arising in
   the level set formulation of the inverse mean curvature flow},
   journal={Adv. Calc. Var.},
   volume={6},
   date={2013},
   pages={123--164},
}
		

\bib{Mawhin-Willem}{book}{
   author={Mawhin, Jean},
   author={Willem, Michel},
   title={Critical point theory and Hamiltonian systems},
   series={Applied Mathematical Sciences},
   volume={74},
   publisher={Springer-Verlag, New York},
   date={1989},
}

\bib{Mercaldo-Segura-Trombetti-2008}{article}{
   author={Mercaldo, A.},
   author={Segura de Le{\'o}n, S.},
   author={Trombetti, C.},
   title={On the behaviour of the solutions to {$p$}-Laplacian equations as
   {$p$} goes to 1},
   journal={Publ. Mat.},
   volume={52},
   date={2008},
   pages={377--411},
}

\bib{Mercaldo-Segura-Trombetti-2009}{article}{
   author={Mercaldo, A.},
   author={Segura de Le{\'o}n, S.},
   author={Trombetti, C.},
   title={On the solutions to 1-Laplacian equation with $L^1$ data},
   journal={J. Funct. Anal.},
   volume={256},
   date={2009},
   pages={2387--2416},
}
		

\bib{Papageorgiou}{article}{
   author={Papageorgiou, Nikolaos S.},
   author={Papageorgiou, Apostolos S.},
   title={Minimization of nonsmooth integral functionals},
   journal={Internat. J. Math. Math. Sci.},
   volume={15},
   date={1992},
   pages={673--679},
}


\bib{Peetre}{article}{
   author={Peetre, Jaak},
   title={Espaces d'interpolation et th\'eor\`eme de Soboleff},
   journal={Ann. Inst. Fourier (Grenoble)},
   volume={16},
   date={1966},
   pages={279--317},
}
\end{biblist}

\end{bibdiv}
\end{document}